\documentclass[12pt]{amsart}
\usepackage[colorlinks=true]{hyperref}
\usepackage[utf8]{inputenc}
\usepackage{amsmath,amsthm,amssymb,url,mathdots,verbatim,psfrag}
\usepackage[margin = 2.5cm]{geometry}
\usepackage{tikz}
\usepackage{ytableau}
\usepackage{color}
\usepackage{MnSymbol}
\usepackage{graphicx}
\usepackage{scalerel}
\usepackage{enumitem}
\usepackage[normalem]{ulem}
\newtheorem{theorem}{Theorem}

\newtheorem{prop}[theorem]{Proposition}

\newtheorem{lemma}[theorem]{Lemma}

\theoremstyle{definition}

\newtheorem{remark}[theorem]{Remark}

\newcommand{\rd}[1]{{\color{red} #1}}

\newcommand{\alert}[1]{\mathbf{\rd{#1}}}

\DeclareMathOperator{\sgn}{sgn}

\newcommand{\nenwarrow}{\nwarrow \!\!\!\!\!\;\!\! \nearrow}
\usepackage{enumitem}

\newcommand{\e}{{\operatorname{E}}}
\newcommand{\W}{{\operatorname{W}}}

\newcommand{\asym}{\mathbf{ASym}}

\newmuskip\pFqmuskip

\newcommand*\pFq[6][8]{%
  \begingroup 
  \pFqmuskip=#1mu\relax
  \mathchardef\normalcomma=\mathcode`,
  \mathcode`\,=\string"8000
  \begingroup\lccode`\~=`\,
  \lowercase{\endgroup\let~}\pFqcomma
  {}_{#2}F_{#3}{\left[\genfrac..{0pt}{}{#4}{#5};#6\right]}%
  \endgroup
}
\newcommand{\pFqcomma}{{\normalcomma}\mskip\pFqmuskip}

\title[Arrowed Gelfand-Tsetlin patterns]{$(-1)$-Enumerations of arrowed Gelfand-Tsetlin patterns}

\author{Ilse Fischer and Florian Schreier-Aigner}
\thanks{The authors acknowledge the financial support from the Austrian Science Foundation FWF, grant P34931}

\begin{document}

\begin{abstract}
Arrowed Gelfand-Tsetlin patterns have recently been introduced to study alternating sign matrices. In this paper, we show that a $(-1)$-enumeration of arrowed Gelfand-Tsetlin patterns can be expressed by a simple product formula. The numbers are up to $2^n$ a one-parameter generalization of 
the numbers $2^{n(n-1)/2} \prod_{j=0}^{n-1} \frac{(4j+2)!}{(n+2j+1)!}$ that appear in recent work of Di Francesco.
A second result concerns the $(-1)$-enumeration of arrowed Gelfand-Tsetlin patterns when excluding double-arrows as decoration in which case we also obtain a simple product formula. We are also able to provide signless interpretations of our results.
 The proofs of the enumeration formulas are based on a recent Littlewood-type identity, which allows us to reduce the problem to the evaluations of two determinants. The evaluations are accomplished by means of the LU-decompositions of the underlying matrices, and an extension of Sister Celine's algorithm as well as creative telescoping to evaluate certain triple sums. In particular, we use implementations of such algorithms by Koutschan, and by Wegschaider and Riese.
\end{abstract}

\maketitle

\section{Introduction} 
Arrowed Gelfand-Tsetlin patterns have appeared in recent work \cite{bounded,nASMDPP} in connection with alternating sign matrices and related Littlewood-type identities.
A \emph{Gelfand-Tsetlin pattern} is a triangular array of integers of the following form  
\begin{equation*} 
\label{triangle} 
\begin{array}{ccccccccccccccccc}
  &   &   &   &   &   &   &   & a_{1,1} &   &   &   &   &   &   &   & \\
  &   &   &   &   &   &   & a_{2,1} &   & a_{2,2} &   &   &   &   &   &   & \\
  &   &   &   &   &   & \dots &   & \dots &   & \dots &   &   &   &   &   & \\
  &   &   &   &   & a_{n-2,1} &   & \dots &   & \dots &   & a_{n-2,n-2} &   &   &   &   & \\
  &   &   &   & a_{n-1,1} &   & a_{n-1,2} &  &   \dots &   & \dots   &  & a_{n-1,n-1}  &   &   &   & \\
  &   &   & a_{n,1} &   & a_{n,2} &   & a_{n,3} &   & \dots &   & \dots &   & a_{n,n} &   &   &
\end{array},
\end{equation*} 
with weak increase along $\nearrow$-diagonals and $\searrow$-diagonals, that is,
$a_{i+1,j} \le a_{i,j} \le a_{i+1,j+1}$ for all $i,j$. \emph{Arrowed Gelfand-Tsetlin patterns} are Gelfand-Tsetlin patterns, where entries may carry
decorations from
$\{\nwarrow, \nearrow, \nenwarrow\}$ such that the following is satisfied: 
\emph{Suppose an entry $a$ is equal to its $\nearrow$-neighbor (resp. $\nwarrow$-neighbor) and $a$ is decorated with either $\nearrow$ or $\nenwarrow$ (resp. $\nwarrow$ or $\nenwarrow$), then the entry right (resp. left) of $a$ in the same row is also equal to $a$ and decorated with $\nwarrow$ or 
$\nenwarrow$ (resp. $\nearrow$ or $\nenwarrow$).}  

In our examples, an entry 
 $a$ that is 
  decorated with $\nwarrow,\nearrow, {\nwarrow \!\!\!\!\!\;\!\! \nearrow}$ appears 
   as $^\nwarrow a$, $a^\nearrow$, $^\nwarrow a^\nearrow$, respectively. Next we provide such an example.
$$
\begin{array}{ccccccccccccccccc}
  &   &   &   &   &   &   &   & ^{\phantom{\nwarrow}}5^\nearrow &   &   &   &   &   &   &   & \\
  &   &   &   &   &   &   & ^\nwarrow 5^{\phantom{\nearrow}} &   & ^\nwarrow7^{\phantom{\nearrow}} &   &   &   &   &   &   & \\
  &   &   &   &   &   & ^{\phantom{\nwarrow}}4^\nearrow &   & \alert{^\nwarrow6^\nearrow} &   & ^{\phantom{\nwarrow}}7^{\phantom{\nearrow}} &   &   &   &   &   & \\
  &   &   &   &   & \alert{^{\phantom{\nwarrow}}4^{\phantom{\nearrow}}} &   & \alert{^{\phantom{\nwarrow}}6 ^\nearrow} &   & \alert{^\nwarrow 6^\nearrow} &   & ^{\phantom{\nwarrow}}7^{\phantom{\nearrow}} &   &   &   &   & \\
  &   &   &   & \alert{^\nwarrow4^\nearrow} &   & \alert{^\nwarrow4^{\phantom{\nearrow}}} &  &   ^{\phantom{\nwarrow}}6^{\phantom{\nearrow}} &   & ^{\phantom{\nwarrow}}7^{\phantom{\nearrow}}   &  & ^\nwarrow9^\nearrow  &   &   &   & \\
  &   &   & ^{\phantom{\nwarrow}}2^\nearrow &   & ^{\phantom{\nwarrow}}4^{\phantom{\nearrow}} &   & ^{\phantom{\nwarrow}}5^\nearrow &   & ^\nwarrow7^{\phantom{\nearrow}} &   & ^\nwarrow8^\nearrow &   & ^{\phantom{\nwarrow}}10^\nearrow &   &   &
\end{array}  
$$

We introduce the \emph{sign} and the \emph{weight} of a given arrowed Gelfand-Tsetlin pattern $A$ as follows:  
Suppose the pattern contains an entry $a$ that is equal to its $\swarrow$-neighbor $b$ as well as to its $\searrow$-neighbor $c$, and $b$ is decorated with $\nearrow$ or $\nenwarrow$ and $c$ is decorated with $\nwarrow$ or $\nenwarrow$, that is the possible arrow patterns for $(b,c)$ are 
$(\nearrow,\nwarrow), (\nearrow,\nenwarrow), (\nearrow,\nwarrow), (\nenwarrow,\nenwarrow)$.   Such a configuration is said to be a \emph{special little triangle}. Then the sign of $A$ is 
$$
\sgn(A):=(-1)^{\# \text{ of special little triangles in } A}.
$$
There are two special little triangles in our example which are highlighted by their bold, red entries, thus the sign is $1$.
Now the weight of $A$ is defined as follows, where $\# \emptyset$ denotes the number of 
entries without decoration.
\begin{equation} 
\label{weight} 
\sgn(A) \cdot t^{\# \emptyset}  u^{\# \nearrow} v^{\# \nwarrow} w^{\# \nenwarrow}   \prod_{i=1}^{n} X_i^{\sum_{j=1}^i a_{i,j}  - \sum_{j=1}^{i-1} a_{i-1,j} + \# \nearrow \text{in row $i$ } - \# \nwarrow \text{in row $i$ }} =: \W(A)
\end{equation} 
In our example, the weight is 
$$
t^{6} u^{6} v^{4} w^{5} X_1^{6} X_2^{5} X_3^{6} X_4^{7} X_5^{6}
X_6^8.
$$

The significance of arrowed Gelfand-Tsetlin patterns partly stems from the following.
If every entry carries a decoration and the bottom row of an arrowed 
Gelfand-Tsetlin pattern is strictly increasing, then the underlying 
undecorated triangular integer array is a \emph{monotone triangle}. 
(A monotone triangle is a
Gelfand-Tsetlin pattern with strictly increasing rows and monotone triangles with bottom row $1,2,\ldots,n$ are in easy bijective correspondence with $n \times n$ \emph{alternating sign matrices}, see \cite{Bre99}.)
Indeed, one can see inductively, from the bottom row to the top row, 
that all rows are strictly increasing as every entry needs to have a decoration and special little triangles are impossible under these circumstances.

Now, if we specialize the generating function of all arrowed Gelfand-Tsetlin patterns with a given bottom row $k_1 < k_2 < \ldots < k_n$ with respect to the above weight by setting 
 $$(X_1,\ldots,X_n)=(1,\ldots,1), t=0, u=v=1, w=-1,$$ 
 then we obtain the number of monotone triangles with bottom row $k_1,\ldots,k_n$. Indeed, setting $t=0$ means that we exclude patterns that have at least one entry without decoration, and, according to the previous paragraph, we are left with monotone triangles.  Given a monotone triangle $M$, what are the eligible decorations, assuming that every entry carries a decoration? Let $d$ be the number of entries in $M$ that are equal to their $\nwarrow$-neighbor (in which case the decoration is 
$\nearrow$), $e$ be the number of entries in $M$ that are equal to their $\nearrow$-neighbor (in which case the decoration is 
$\nwarrow$) and  $f$ be the number of the remaining entries in $M$ (such entries may be decorated with any of the three options). Then the sum of all weights of arrowed Gelfand-Tsetlin patterns that have $M$ as underlying monotone triangle is --- when setting $(X_1,\ldots,X_n)=(1,\ldots,1)$ --- equal to $u^{d} v^{e} (u+v+w)^f$.
After specializing $u=v=1$ and $w=-1$, the weight is obviously $1$.

We come back to arrowed Gelfand-Tsetlin patterns that have also entries without decorations and present the results of this paper. We compute for fixed $n$ and $m$, the generating function of arrowed Gelfand-Tsetlin patterns with $n$ rows and all entries are (possibly decorated) non-negative integers no greater than $m$ such that the bottom row is strictly increasing, for the following specialization
$$(X_1,\ldots,X_n)=(1,\ldots,1), u=v=t=1, w=-1.$$
In other words, we deal with the $(-1)$-enumeration of arrowed Gelfand-Tsetlin patterns with respect to the weight 
$$
\sgn(A) (-1)^{\# \nenwarrow} = (-1)^{\# \text{special little triangles} + \# \nenwarrow}.
$$

\begin{theorem} 
\label{1}
The specialization of the generating function of arrowed Gelfand-Tsetlin patterns with $n$ rows and strictly increasing non-negative bottom row where the entries are bounded by $m$ at 
$(X_1,\ldots,X_n)=(1,\ldots,1)$, $u=v=t=1$ and $w=-1$ is
equal to 
$$
2^n \prod_{i=1}^n \frac{(m-n+3i+1)_{i-1} (m-n+i+1)_i}{\left(\frac{m-n+i+2}{2} \right)_{i-1} (i)_i}, 
$$
using Pochhammer notation $(a)_n=a(a+1) \dots (a+n-1)$.
\end{theorem}

When setting  $m=n-1$ in this theorem, we obtain up to a power of $2^n$ the following numbers.
\begin{equation}
\label{difran}
1,4,60,3328,678912,\ldots = 2^{n(n-1)/2} \prod_{j=0}^{n-1} \frac{(4j+2)!}{(n+2j+1)!}=F_n
\end{equation}
Indeed, when setting $m=n-1$ in the formula in the theorem, we obtain 
$$
2^n \prod_{i=1}^n \frac{((n-1)-n+3i+1)_{i-1} ((n-1)-n+i+1)_i}{\left(\frac{(n-1)-n+i+2}{2} \right)_{i-1} (i)_i} = 
2^n \prod_{i=1}^n \frac{(3i)_{i-1}}{\left(\frac{i+1}{2} \right)_{i-1}}=:F'_n,   
$$
so that 
$\frac{F'_{n}}{F'_{n-1}}=2 \frac{(3n)_{n-1}}{\left( \frac{n+1}{2} \right)_{n-1}}$, 
while 
$$
\frac{F_{n}}{F_{n-1}} = 2^{n-1} \frac{(4n-2)!}{(3n-1)!} \prod_{j=0}^{n-2} \frac{1}{n+2j+1},  
$$
which is up to the factor $2$ equal to $\frac{F'_{n+1}}{F'_{n}}$, and 
$F'_{n} = 2^n F_n$ now follows from $F_n=1$ and $F'_n=2$. 
The significance of this is that the numbers in \eqref{difran} have also appeared in recent work of Di Francesco \cite{twenty}. Di Francesco conjectured that these numbers count $20$-vertex configurations on the quadrangle as well as domino tilings of the Aztec triangle (for the precise definition we refer to Di Francesco's paper). Meanwhile the conjecture has been proved in a recent article by Koutschan, Krattenthaler and Schlosser \cite{KoutschanKrattenthalerSchlosser24plus}.
 Setting $m=n-1$ for the combinatorial objects in Theorem~\ref{1} means that the bottom row of the arrowed Gelfand-Tsetlin pattern is prescribed as $0,1,\ldots,n-1$. This concerns an earlier unpublished conjecture from 2018 of the second author.

Our second theorem concerns the $(-1)$-enumeration of arrowed Gelfand-Tsetlin patterns, where we exclude $\nenwarrow$ as decoration.

\begin{theorem} 
\label{3}
The specialization of the generating function of arrowed Gelfand-Tsetlin patterns with $n$ rows and strictly increasing non-negative bottom row where the entries are bounded by $m$ at 
$(X_1,\ldots,X_n)=(1,\ldots,1)$, $u=v=t=1$ and $w=0$ is
equal to 
$$
3^{\binom{n+1}{2}} \prod_{i=1}^n \frac{(2n+m+2-3i)_i}{(i)_i}.
$$ 
\end{theorem} 

Note that in this case, the weight is simply 
$$\sgn(A)=(-1)^{\# \text{ of special little triangles in } A}.$$

\subsection*{Outline of the paper} In Section~\ref{signless}, we provide signless combinatorial interpretations of the numbers in Theorems~\ref{1} and \ref{3} in terms of weighted counts of some (possibly decorated) Gelfand-Tsetlin patterns with positive weights. In Section~\ref{prel}, we summarize necessary results from other papers. In Section~\ref{sec:odd}, we prove Theorem~\ref{1}, and, in Section~\ref{sec:3}, we prove Theorem~\ref{3}.

\section{Signless versions of the objects} 
\label{signless} 

In this section, we show that our $(-1)$-enumerations can also be expressed as weighted counts of certain (possibly decorated) Gelfand-Tsetlin patterns with respect to \emph{positive} weights. The first proposition concerns the $(-1)$-enumeration in Theorem~\ref{1}.

\begin{prop}
The $(-1)$-enumeration of arrowed Gelfand-Tsetlin patterns with given bottom row $k_1 \le \ldots \le k_n$ with respect to the sign 
$$
(-1)^{\# \text{special little triangles} + \# \nenwarrow}
$$
is the weighted enumeration of Gelfand-Tsetlin patterns with bottom row $k_1 \le \ldots \le k_n$ such that only 
the bottom row can contain three equal entries and where the weight is 
$2^{r}$ such that 
$r$ is the number of entries in the Gelfand-Tsetlin pattern that are not equal to both their $\nwarrow$-neighbor  \emph{and} their $\nearrow$-neighbor. 
\end{prop}  

Note that in this proposition, we are using a notion of sign of arrowed Gelfand-Tsetlin pattern that is different from $\sgn(A)$.

\begin{proof}
We start by giving a \emph{sign-reversing involution on arrowed Gelfand-Tsetlin patterns that have a row with at least four equal entries.} Choose the topmost  and leftmost occurrence of four equal entries in a row, and denote them by $a_1,a_2,a_3,a_4$, from left to right. They are all adjacent and the $\nwarrow$-neighbor of $a_1$ is strictly less than $a_1$ and the $\nearrow$-neighbor 
of $a_4$ is strictly greater than $a_4$ because otherwise $a_1,a_2,a_3,a_4$ were not a topmost occurrence of four equal entries. If $a_1,a_2$ are part of 
the same special little triangle, then a sign-reversing involution is to change the decoration of $a_1$ from $\nearrow$ to $\nenwarrow$, or vice versa. Similar for $a_3$ and $a_4$.
Thus we can assume that $a_1$ is either not decorated or decorated with $\nwarrow$, 
and that $a_4$ is either not decorated or decorated with $\nearrow$. Therefore, $a_2$ is either not decorated or decorated with $\nearrow$, while $a_3$ is either not decorated or decorated with $\nwarrow$. More precisely, $a_2$ is decorated with $\nearrow$ if and only if $a_3$ is decorated with $\nwarrow$. Therefore, we have a sign-reversing involution by changing from no decoration on $a_2$ and $a_3$ to decorating $a_2$ with $\nearrow$ and $a_3$ with $\nwarrow$, since $a_2$
and $a_3$ are part of the same special little triangle in the second case but not in the first case. From now on, we can assume that there are at most three equal entries in a row.

Next assume that \emph{there are three equal entries $a_1,a_2,a_3$ in a row, from left to right, that are not in the bottom row, and we will also construct a sign-reversing involution in this case.} We assume that $a_1,a_2,a_3$ is the topmost and leftmost occurrence of three equal entries in a row.
Let $b_1,b_2$ be the entries in the row below such that $a_1$ is the $\nwarrow$-neighbor of $b_1$, and $b_2$ is the right neighbor of $b_1$. 
We show that the left neighbor of $b_1$ is different from $b_1$ and that 
the right neighbor of $b_2$ is different from $b_2$. Since there are no four equal entries in a row it cannot happen that we have equality on both sides.
Assume without of generality that the right neighbor $b_3$ of $b_2$ is equal to $b_2$. We define a sign-reversing involution of arrowed Gelfand-Tsetlin pattern such that $b_2$ and $b_3$ are contained in the same special little triangle by changing the decoration of $b_3$ from $\nwarrow$ to 
$\nenwarrow$, and vice versa. This is possible as the $\nearrow$-neighbor $a_4$ of $b_3$ is different from $b_3$ because otherwise there would be four equal entries in a row, namely $a_1,a_2,a_3,a_4$. This implies that $b_2$ is either not decorated or decorated with 
$\nwarrow$. Now $b_1$ is either not decorated or decorated with $\nearrow$ as its left neighbor is different from $b_1$, since there are no four equal entries in a row and $b_1$ is equal to its $\nwarrow$-neighbor. Moreover, $b_1$ is decorated with $\nearrow$ if and only of $b_2$ is decorated with $\nwarrow$, and we obtain a sign reversing 
involution by changing this decoration to no decoration of $b_1$ and $b_2$, and vice versa.
It follows that we can restrict to Gelfand-Tsetlin patterns where three equal entries may only appear in the bottom row.

Now suppose that we have a special little triangle and let $a_1,a_2$ be its two bottom entries. Then either the $\nwarrow$-neighbor of $a_1$ is different from $a_1$, or the $\nearrow$-neighbor of $a_2$ is different from $a_2$ because otherwise we would have three equal entries in the row above the row of $a_1,a_2$. Without loss of generality, we may assume that 
$a_1$ is different from its $\nwarrow$-neighbor and we have a sign-reversing involution 
by changing the decoration of $a_1$ from $\nwarrow$ to $\nenwarrow$, and vice versa. \emph{From now on we can assume that there are no special little triangles.}

It follows that the specialization gives the signed enumeration of decorated Gelfand-Tsetlin patterns such that only the bottom row can contain 
three equal entries, and  the entries may carry decorations from
$\{\nwarrow, \nearrow, \nenwarrow\}$ such that if the entry $a$ is equal to its $\nearrow$-neighor (resp. $\nwarrow$-neighbor), $a$ must not carry the decoration $\nearrow$ or $\nenwarrow$ (resp.
$\nwarrow$ or $\nenwarrow$). The sign is $(-1)^{\# \nenwarrow}$.

An entry that is equal to its $\nwarrow$-neighbor and its $\nearrow$-neighbor cannot be decorated (since there is no special little triangle), thus we may think of such an entry to be weighted by $1$ (this corresponds to no decoration). An entry that is equal to its $\nwarrow$-neighbor but not to its $\nearrow$-neighbor (if it exists) can only be decorated with $\nearrow$
or has no decoration, thus it is weighted by $2$, likewise an entry that is equal to its $\nearrow$-neighbor but not to its $\nwarrow$-neighbor (if it exists) can only be decorated with $\nwarrow$ or has no decoration, thus it is also weighted by $2$. For all other entries, all decorations are possible, and, since $t+u+v+w=2$ under our specialization, they are also weighted by $2$.
\end{proof} 

Note that in our application in Theorem~\ref{1}, we have $k_1 < k_2 < \ldots < k_n$ and thus we consider Gelfand-Tsetlin patterns that have at most two equal entries in a row. 

\medskip

We turn to the $(-1)$-enumeration from Theorem~\ref{3}.

\begin{prop} 
The $(-1)$-enumeration of arrowed Gelfand-Tsetlin patterns where no entry is decorated with $\nenwarrow$
with given bottom row $k_1 \le \ldots \le k_n$ is 
the generating function of Gelfand-Tsetlin patterns with bottom row $k_1 \le \ldots \le k_n$
such that each entry appears at most twice in each row and entries may carry decorations from 
$\{\nwarrow,\nearrow\}$ so that the following is satisfied. 
\begin{itemize}
\item If the entry $a$ is equal to its $\nearrow$-neighbor (resp. $\nwarrow$-neighbor), then $a$ must not carry the decoration $\nearrow$  (resp. $\nwarrow$), and
\item if $a,b$ are two equal entries in the same row, then at least one of $a$ and $b$ has to be decorated.
\end{itemize} 
\end{prop}

\begin{proof} 
Letting $w=0$ means that we exclude the decoration $\nenwarrow$. Let $a_1, a_2, a_3$ be the topmost and leftmost occurrence of three equal entries in a row. It follows that the $\nwarrow$-neighbor of $a_1$ cannot be equal to $a_1$ and the $\nearrow$-neighbor of $a_3$ cannot be equal to $a_3$ because otherwise $a_1, a_2, a_3$ would not be the topmost occurrence of three equal entries. There are the following eight possible decorations of $a_1,a_2,a_3$:
$$
(\emptyset,\emptyset,\emptyset),(\emptyset,\emptyset,\nearrow),
(\emptyset,\nearrow,\nwarrow),
(\nwarrow,\emptyset,\emptyset),(\nwarrow,\emptyset,\nearrow),
(\nwarrow,\nearrow,\nwarrow),
(\nearrow,\nwarrow,\emptyset),(\nearrow,\nwarrow,\nearrow)
$$
and the signs that they contribute are $+,+,-,+,+,-,-,-$ in that order, and we obtain a sign-reversing involution. Therefore we can assume that there are at most two occurrences of the same integer in a row.

Another sign-reversing involution is obtained in the set of the remaining arrowed Gelfand-Tsetlin patterns that have either a special little triangle or two equal entries in the same row that are both not decorated as follows: Consider the topmost and leftmost occurrence of either of these local configurations and transform one configuration into the other.
\end{proof}

\section{Preliminaries for the proof}
\label{prel} 

In this section, we summarize results from other papers that will be used for the proof. For reasons of readability, we omit writing the summation bounds whenever possible, that is, when the sum can be extended to $\mathbb{Z}$ since the summand 
vanishes for small and large parameters. In particular, all summations are actually finite and we write $\sum\limits_i$ when summing over all $i \in \mathbb{Z}$.
\bigskip

It follows from \cite[Theorem~3.2]{nASMDPP} and manipulations that are analogous to those in \cite[Eq.\ (7.2)]{nASMDPP} that the generating function of arrowed Gelfand-Tsetlin patterns with $n$ rows and entries at most $m$ is given by 
$$
\frac{1}{\prod\limits_{1 \le i < j \le n} (X_j-X_i)} \asym_{X_1,\ldots,X_n} \left[  \prod_{1 \le i \le j \le n} (1+w X_i+X_j + X_i X_j) 
\sum_{0 \le k_1 < k_2 < \ldots < k_n \le m} X_1^{k_1-1} X_2^{k_2-1} \cdots X_n^{k_n-1} \right],  
$$
where the \emph{antisymmetrizer} $\asym$ is defined as 
$$
\asym_{X_1,\ldots,X_n} f(X_1,\ldots,X_n) = \sum_{\sigma \in {\mathcal S}_n} \sgn \sigma 
\cdot f(X_{\sigma(1)},\ldots,X_{\sigma(n)}),
$$
see also \cite[Eq.\ (3.5)]{bounded}.
Then, in \cite[Cor.\ 1.2]{bounded}, it was shown that this is equal to 
\begin{multline}
\label{bialternant} 
\prod_{i=1}^n (X_i^{-1} + 1 + w + X_i) \\ \times
\frac{\det_{1 \le i, j \le n} \left( X_i^{j-1}  (1+X_i)^{j-1} (1+ w X_i)^{n-j}  - X_i^{m+2n-j} (1+X_i^{-1})^{j-1}  (1+w X_i^{-1})^{n-j} \right)}{\prod\limits_{i=1}^n (1-X_i) \prod\limits_{1 \le i < j \le n} (1-X_i X_j)(X_j-X_i)}.
\end{multline} 
In that paper, we also gave a formula in terms of complete homogeneous symmetric polynomials $h_k$. In order to formulate it, we need the $h_k$ to be extended to negative $k$ as follows. We define $h_k(X_1,\ldots,X_n)=0$ for $k=-1,-2,\ldots,-n+1$ and 
$$
h_k(X_1,\ldots,X_n) = (-1)^{n+1} X_1^{-1} \ldots X_n^{-1} h_{-k-n}(X_1^{-1},\ldots,X_n^{-1})
$$
for $k \le -n$. A consequence of this definition is that the latter relation is true for any $k$.
For the formula in terms of the complete symmetric functions, we need to distinguish between the cases $m$ is odd and $m$ is even. At the end of this section, we argue why it suffices to consider the case that $m$ is odd for our results. Now, for $m=2l+1$, \eqref{bialternant} is equal to
\begin{multline} 
\label{odd} 
(-1)^{\binom{n+1}{2}} \prod_{i=1}^n X_i^l (X_i^{-1}+1+w+X_i)(1+X_i) \\
\times \det_{1 \le i, j \le n} \left( \sum_{p,q} \binom{j-1}{p} \binom{n-j}{q} w^{n-j-q} 
h_{p-q-l-1-i}(X_1,X_1^{-1},\ldots,X_i,X_i^{-1}) \right), 
\end{multline} 
(see \cite[(4.9)]{bounded}). The sum that constitutes the entry underlying the determinant has natural bounds: this is due to the binomial coefficients as they vanish if $p$ or $q$ are negative, or $p > j-1$, or $q > n-j$.  

We aim at specializing $(X_1,\ldots,X_n)=(1,\ldots,1)$. For this purpose, observe that 
$$
\left. h_k(X_1,\ldots,X_n) \right|_{(X_1,\ldots,X_n)=(1,\ldots,1)} = \binom{k+n-1}{n-1}
$$
for $k \ge 0$, and that the identity $\binom{p}{q} = (-1)^{q} \binom{q-p-1}{q}$ for the binomial coefficient implies that the formula is also true for $k<0$, 
where we use $\binom{p}{q} = \frac{p(p-1) \dots (p-q+1)}{q!}$.
For the case $m=2l+1$, we obtain the following specialization
\begin{equation}
\label{odddet}
(3+w)^n 2^n \det_{1 \le i, j \le n} \left( \sum_{p,q} w^{n-j-q} (-1)^j \binom{j-1}{p} \binom{n-j}{q}  
\binom{p-q-l+i-2}{2i-1} \right).
\end{equation}

\subsection*{The case $m$ is even}
In the following we argue that it suffices to prove Theorems~\ref{1} and \ref{3} for odd values of $m$. 

Theorem~3.2 from \cite{nASMDPP} implies that the generating function of arrowed Gelfand-Tsetlin patterns with bottom row $k_1,\ldots,k_n$ 
is 
\begin{equation} 
\label{poly} 
\prod_{i=1}^n (t + u X_i + v X_i^{-1} + w) \prod_{1 \le p < q \le n} (t + u \e_{k_p} + v \e_{k_q}^{-1} + w \e_{k_p} \e_{k_q}^{-1}) 
s_{(k_n,k_{n-1},\ldots,k_1)}(X_1,\ldots,X_n), 
\end{equation} 
where $s_{(k_n,k_{n-1},\ldots,k_1)}(X_1,\ldots,X_n)$ denotes the extended Schur polynomial, that is 
$$
s_{(k_n,\ldots,k_1)}(X_1,\ldots,X_n) = \frac{\det_{1 \le i,j \le n} (X_i^{k_j+j-1})}{\prod_{1 \le i < j \le n} (X_j - X_i)}
$$
and $\e_x$ denotes the forward shift, that is $\e_x p(x) = p(x+1)$. The formula has to be interpreted as follows: First we treat 
the $k_i$'s as variables and apply the polynomial in the shift operators to the extended Schur polynomial. We obtain a linear combination of Schur polynomials and only then specialize the $k_i$'s. 

Setting $X_i=1$ in $s_{(k_n,k_{n-1},\ldots,k_1)}(X_1,\ldots,X_n)$ 
gives $\prod_{1 \le i < j \le n} \frac{k_j-k_i+j-i}{j-i}$. For $0 \le k_1 \le k_2 \le \ldots \le k_n$, this follows from \cite[p. 375, Eq. (7.105)]{Sta99}. Now, to prove it for arbitrary $k_i$, observe that if there are $i, j$ with $1 \le i \not= j \le n$, $k_i+i=k_j+j$ then both expressions vanish. Otherwise, there is a permutation $\sigma$ such that 
$$
k_{\sigma(1)} + \sigma(1)   < k_{\sigma(2)} + \sigma(2)   < \ldots < k_{\sigma(n)} + \sigma(n) 
$$
and an integer $p$ with $k_1+\sigma(1)-1+p \ge 0$, and we have  
\begin{multline*} 
s_{(k_n,\ldots,k_1)}(X_1,\ldots,X_n) = \sgn(\sigma) \prod_{i=1}^n X_i^{-p} \\
\times s_{(k_{\sigma(n)}+\sigma(n)-n+p,\ldots,k_{\sigma(2)}+\sigma(2)-2+p,k_{\sigma(1)}+\sigma(1)-1+p)}(X_1,\ldots,X_n), 
\end{multline*} 
which reduces it to the previous case. It follows that specializing at $(X_1,\ldots,X_n)=(1,\ldots,1)$ and $t=u=v=1$ 
 in \eqref{poly} gives a polynomial in $k_i$ over $\mathbb{Q}[w]$. Summing over all $k_1,\ldots,k_n$ with $0 \le k_1 < k_2 < \ldots < k_n \le m$ gives then a polynomial in $m$ over $\mathbb{Q}[w]$. This follows inductively from the fact that for a polynomial $p(x)$ we have that $\sum_{x=0}^y p(x)$ is a polynomial function in $y$.
 
Setting $l=\frac{m-1}{2}$ in \eqref{odddet}, we obtain a polynomial in $m$ for fixed $n$, and we know that it agrees with the weighted count of arrowed monotone triangle with entries no greater than $m$ and with respect to the weight 
$$(-1)^{\# \text{special little triangles}} w^{\# \nenwarrow}.$$
Therefore it must agree on all integer values of $m$.
The formulas in Theorems~\ref{1} and \ref{3} are obviously rational functions in $m$. In the following two sections we will show that 
these rational functions agree with \eqref{odddet}, when replacing $l=\frac{m-1}{2}$ in the determinant, if $m$ is odd in the cases $w=-1$ and $w=0$, respectively. Since 
these rational functions agree with polynomials for an infinite number of values (the odd integers), they must be polynomials themselves and thus must also agree with \eqref{odddet} for all values of $m$.

\section{Proof of Theorem~\ref{1}} 
\label{sec:odd} 

\subsection{LU-decomposition} 
\label{sec:LU-decomposition}
We start by modifying the entry of the matrix in \eqref{odddet} further when $w=-1$:
$$
\sum_{p,q} (-1)^{n+q} \binom{j-1}{p} \binom{n-j}{q}  
\binom{p-q-l+i-2}{2i-1} = \sum_{p,q} (-1)^{n+q+1} \binom{j-1}{p} \binom{n-j}{q}  
\binom{l-p+q+i}{2i-1}, 
$$
where we have used $\binom{p}{q} = (-1)^{q} \binom{q-p-1}{q}$ for the last binomial coefficient.
Since this is a polynomial in $l$ and we can therefore assume that $l$ is big enough, we may apply the symmetry of the binomial coefficient.
We also use the Chu-Vandermonde summation to obtain the following. 
\begin{multline*} 
\sum_{p,q} (-1)^{n+q+1} \binom{j-1}{p} \binom{n-j}{q}  
\binom{l-p+q+i}{l-p+q-i+1} \\
= \sum_{p,q} (-1)^{i+l+n+p} \binom{j-1}{p} \binom{n-j}{n-j-q}  
\binom{-2i}{l-p+q-i+1} \\ 
=
\sum_{p} (-1)^{i+l+n+p} \binom{j-1}{p} \binom{n-2i-j}{n-j+l-p-i+1} \\
= (-1)^{j+1} \sum_{p}  \binom{j-1}{p} \binom{l-p+i}{n-j+l-p-i+1} = 
(-1)^{j+1} \sum_{p}  \binom{j-1}{p} \binom{l-p+i}{2i+j-n-1}
\end{multline*}
Since the last expression is again a polynomial in $l$ and since it agrees with the original one for large $l$, it must agree for all $l$. 
Now, since this is the entry of a matrix we want to compute the determinant of, we can also replace $j$ by $n+1-j$ by taking into account the sign 
$(-1)^{\binom{n}{2}}$. In summary, we want to compute 
\begin{equation} 
\label{odddet1} 
2^{2n} \det_{1 \le i, j \le n} \left( \sum_{p} \binom{n-j}{p} \binom{l-p+i}{2i-j} \right).
\end{equation}

We define 
$$
a_{i,j} = \sum_{p} \binom{n-j}{p} \binom{l-p+i}{2i-j}
$$
and
\begin{equation}
x_{i,j} = 
\begin{cases} 
\begin{aligned} 
& (-1)^{i+1} \frac{(j)_j}{(2l-n+3j+2)_{j-1} (2l-n+i+2)_j}  \\ 
& \times
\sum_{t} 2^{2i-4t-n} \frac{(i-j-2t+1)_{2t} (i-2j+1)_{j-1-t} (l-n/2+j/2+t+3/2)_{i-2t-1}}{(1)_t (1)_{i-2t-1}}  
\end{aligned} & i \le j \\
0 & \text{otherwise} 
\end{cases}, 
\end{equation}
using the extended Pochhammer  notation
$$(a)_k = \begin{cases} a(a+1) \dots (a+k-1), & k > 0 \\ 
                                       1, & k=0 \\
                                      \frac{1}{(a-1)(a-2) \cdots (a+k)}, & k<0 
                                      \end{cases}.$$
The remainder of the section is devoted to the proof of the following lemma.

\begin{lemma}
\label{prop:w1 matrix product}
Let $n$ be a positive integer, then the product $(a_{i,j})_{1 \le i,j \le n} \cdot (x_{i,j})_{1 \le i,j \le n}$ is a lower 
triangular matrix with $1$'s on the main diagonal, that is 
\begin{equation}
\label{eq:w1 matrix identity}
\sum_{k=1}^n a_{i,k}  x_{k,j} = \begin{cases}
1 \qquad &i=j,\\
0 & i<j.
\end{cases}
\end{equation}
\end{lemma}
The lemma then implies 
$$
\det_{1 \le i, j \le n} \left( \sum_p \binom{n-j}{p} \binom{l-p+i}{2i-j} \right) = \left[ \det_{1 \le i, j \le n} (x_{i,j}) \right]^{-1}, 
$$
and it is easy to see that $2^{2n}\left[ \det_{1 \le i, j \le n} (x_{i,j}) \right]^{-1}$ is the expression in Theorem~\ref{1} when setting 
$m=2l+1$ there. 

The following sketches the steps of the proof of Lemma~\ref{prop:w1 matrix product}. 
\begin{itemize}
\item In Section~\ref{sec:applying hyp geom}, hypergeometric transformation formulas are used to transform \eqref{eq:w1 matrix identity} into the triple sum identity in \eqref{triple}.
\item In Section~\ref{sec:j large}, we show \eqref{triple} for $j\ge 2i$. This will serve as the  base case of an induction appearing at the end of Section~\ref{sec:proof of triple sum}.
\item In Section~\ref{sec:auxiliary double sum}, we evaluate an auxiliary double sum, which is also needed in the final section.
\item In Section~\ref{sec:proof of triple sum}, we first establish a recursion for the summand of the triple sum in \eqref{triple}, which  leads to a recursion for the sum itself. The recursion can then be simplified using the evaluation in Section~\ref{sec:auxiliary double sum}. We then use induction to prove the assertion.
\end{itemize}

Also note that below we apply a hypergeometric transformation formula to transform the entry $x_{i,j}$ into another hypergeometric sum and we could have simply defined $x_{i,j}$ as the outcome of the transformation. However, $x_{i,j}$, as defined here, is the result of a guessing procedure and it might be of some value to keep this expression.

\subsection{Application of some hypergeometric transformation formulas} 
\label{sec:applying hyp geom}
We start by modifying the entry $\sum_p \binom{n-j}{p} \binom{l-p+i}{2i-j}$ by using a transformation for hypergeometric series. We use the standard 
hypergeometric notation, namely 
$$\pFq{p}{q}{a_1,a_2,\ldots,a_p}{b_1,b_2,\ldots,b_q}{z} = \sum_{k} \frac{(a_1)_k (a_2)_k \dots (a_p)_k}{(b_1)_k (b_2)_k \dots (b_q)_k} \frac{z^k}{k!}.$$
For the manipulations of hypergeometric series we have used Krattenthaler's Mathematica package HYP \cite{hyp}.

Now it can be checked that 
\begin{equation}
\label{eq:w1 hyp manipulation1}
\sum_{p} \binom{n-j}{p} \binom{l-p+i}{2i-j} = \frac{(1-i+j+l)_{2i-j}}{(1)_{2i-j}} \pFq{2}{1}{i-j-l,j-n}{-i-l}{-1}. 
\end{equation}
Using the following transformation formula 
\cite[(1.8.10)]{slater} 
$$
\pFq{2}{1}{a,-n}{c}{z} = (1-z)^{n} \frac{(a)_n}{(c)_n} \pFq{2}{1}{-n,c-a}{1-a-n}{\frac{1}{1-z}},
$$
where $n$ is a non-negative integer, \eqref{eq:w1 hyp manipulation1} becomes 
\begin{multline*} 
2^{-j+n} \frac{(1-i+j+l)_{2i-j}(i-j-l)_{n-j}}{(1)_{2i-j}(-i-l)_{n-j}}\pFq{2}{1}{-2i+j,j-n}{1-i+2j+l-n}{\frac{1}{2}} \\ =
2^{-j+n} \frac{(1-i+2j+l-n)_{2i-j}}{(1)_{2i-j}}\pFq{2}{1}{-2i+j,j-n}{1-i+2j+l-n}{\frac{1}{2}} \\ 
= \sum_{s} (-1)^s 2^{-j+n-s} \frac{(j-n)_s(1-i+2j+l-n+s)_{2i-j-s}}{(1)_s (1)_{2i-j-s}},
\end{multline*} 
where $n-j$ serves as the non-negative integer and Pochhammer symbols have been combined in the last two steps. By the above, \eqref{eq:w1 matrix identity} is equivalent to
\begin{multline*} 
\sum_{s,t \atop 1 \le k \le j} \frac{(1-i+2k+l-n+s)_{2i-k-s} (l-n/2+j/2+t+3/2)_{k-2t-1}}{(2l-n+3j+2)_{j-1} (2l-n+k+2)_j} \\
\times (-1)^{k+s+1} 2^{k-s-4t} \frac{(k-n)_s (j)_j (k-j-2t+1)_{2t} (k-2j+1)_{j-1-t} }{(1)_s (1)_{2i-k-s}  (1)_t (1)_{k-2t-1} }  = 
\begin{cases} 1, & i=j, \\ 0, & i<j. \end{cases} 
\end{multline*} 
We multiply this identity with $\frac{(2l-n+3)_{2j-1} (2l-n+3j+2)_{j-1}}{(j)_j}$ and replace $l$ by $(n-r-3)/2$. 
\begin{multline*}
 \sum_{s,t \atop 1 \le k \le j} (-r)_{k-1} (-1+j+k-r)_{j-k} \left(-1/2 - i + 2k - n/2 - r/2 + s \right)_{2i-k-s} (j/2-r/2+t)_{k-2t-1} \\
\times (-1)^{k+s+1} 2^{k-s-4t} \frac{(k-n)_s (k-j-2t+1)_{2t} (k-2j+1)_{j-1-t} }{(1)_s (1)_{2i-k-s}  (1)_t (1)_{k-2t-1} }  \\ 
= 
\begin{cases} \frac{(-r)_{2j-1} (-1+3j-r)_{j-1}}{(j)_j}, & i=j, \\ 0, & i<j. \end{cases} 
\end{multline*} 
We aim at applying a transformation rule for the sum over $t$. For this purpose, it suffices to consider the part of the summand that depends only on 
$t$. 
\begin{multline} 
\label{eq:sum over t}
\sum_{t} 2^{-4t} \frac{(j/2-r/2+t)_{k-2t-1} (k-j-2t+1)_{2t} (k-2j+1)_{j-1-t}}{(1)_t (1)_{k-2t-1}} \\
= \frac{(1-2j+k)_{j-1} \left( \frac{j}{2}-\frac{r}{2} \right)_{k-1}}{(1)_{k-1}} \pFq{4}{3}{\frac{1}{2}-\frac{k}{2},1-\frac{k}{2},\frac{j}{2}-\frac{k}{2},\frac{1}{2}+\frac{j}{2}-\frac{k}{2}}{1+j-k,2-\frac{j}{2}-k+\frac{r}{2},\frac{j}{2}-\frac{r}{2}}{1}
\end{multline} 
We apply 
\begin{multline*} 
\pFq{4}{3}{a,b,c,-n}{e,f,1+a+b+c-e-f-n}{1} \\
= \frac{(-a+e)_n (-a+f)_n}{(e)_n (f)_n}  
 \pFq{4}{3}{-n,a,1+a+c-e-f-n,1+a+b-e-f-n}{1+a+b+c-e-f-n,1+a-e-n,1+a-f-n}{1},
\end{multline*} 
where $n$ is a non-negative integer, see \cite[(4.3.5.1)]{slater}, 
and obtain
\begin{multline*} 
\frac{(1-2j+k)_{j-1} \left( \frac{j}{2}-\frac{r}{2} \right)_{k-1}}{(1)_{k-1}} 
\pFq{4}{3}{\frac{1}{2}-\frac{k}{2},1-\frac{k}{2},-\frac{1}{2}+\frac{k}{2}-\frac{r}{2},\frac{k}{2}-\frac{r}{2}}
{\frac{j}{2}-\frac{r}{2},\frac{1}{2}+\frac{j}{2}-\frac{r}{2},\frac{3}{2}-j}{1} \\ \times 
\frac{(1/2+j-k/2)_{k/2-1} (3/2-j/2-k/2+r/2)_{k/2-1}}{(1+j-k)_{k/2-1}(2-j/2-k+r/2)_{k/2-1}}
\end{multline*}
if $k$ is even and 
\begin{multline*} 
\frac{(1-2j+k)_{j-1} \left( \frac{j}{2}-\frac{r}{2} \right)_{k-1}}{(1)_{k-1}} 
\pFq{4}{3}{\frac{1}{2}-\frac{k}{2},1-\frac{k}{2},-\frac{1}{2}+\frac{k}{2}-\frac{r}{2},\frac{k}{2}-\frac{r}{2}}
{\frac{j}{2}-\frac{r}{2},\frac{1}{2}+\frac{j}{2}-\frac{r}{2},\frac{3}{2}-j}{1} \\ \times 
\frac{(j-k/2)_{k/2-1/2} (1-j/2-k/2+r/2)_{k/2-1/2}}{(1+j-k)_{k/2-1/2}(2-j/2-k+r/2)_{k/2-1/2}}
\end{multline*}
if $k$ is odd. Combining both cases, we see that \eqref{eq:sum over t} is equal to
$$
\sum_{t} (-1)^{j+k+t} 2^{-2k+2} \frac{(1-k)_{2t} (-1+k-r)_{2t} (j-t)_{j-k+t}}{(1)_{k-1} (1)_t (1-2j+k)_{1-k+2t} (-1+j+k-r)_{1-k+2t}}.
$$ 
It follows that, in total, \eqref{eq:w1 matrix identity} is equivalent to 
\begin{multline} 
\label{triple}
\sum_{s,t \atop 1 \le k \le j} (-1)^{1+k+s+t} 2^{2-k-s} \\
\times \frac{(k-n)_s (-r)_{k-1+2t} (-1/2-i+2k-n/2-r/2+s)_{2i-k-s} (2-2j+r)_{j-1-2t} (j-t)_{j-t-1}}{(1)_{2i-k-s} (1)_s (1)_{k-1-2t} (1)_t} \\ 
 = \begin{cases} \frac{(-r)_{2j-1} (-1+3j-r)_{j-1}}{(j)_j}, & 0<i=j, \\ 0, & 0<i<j \end{cases} 
\end{multline}

To find the proofs of this triple sum, we were also assisted by Wegschaider's and Riese's MultiSum package \cite{multisum}, which is available at \url{https://www3.risc.jku.at/research/combinat/software/ergosum/RISC/MultiSum.html}. It was useful to find recursions for serveral summands. All recursions that were found by this package are proved in the paper as they are not complicated. 

\subsection{The case $j \ge 2i$}
\label{sec:j large}
We first consider the case $j \ge 2i$, which serves as the base case for the induction at the end of Section~\ref{sec:proof of triple sum}. In this case, we can extend the sum over all $k \in \mathbb{Z}$. Indeed, for $k<1$ either the factor $1/(1)_{k-1-2t}$ vanishes for $t\geq 0$ or $1/(1)_{t}$ vanishes for $t<0$; for $k>j \ge 2i$ either $1/(1)_{2i-k-s}$ vanishes for $s \ge 0$ or $1/(1)_s$ vanishes for $s<0$.

It turns out that in this case already the double sum over $k,s$ vanishes. Neglecting all factors in \eqref{triple} that are independent of $k$ and $s$, it suffices to show 
\begin{equation} 
\label{sum1}
\sum_{k,s} (-1)^{k+s} 2^{-k-s} \frac{(k-n)_s (-r)_{k-1+2t} (-1/2-i+2k-n/2-r/2+s)_{2i-k-s}}{(1)_{2i-k-s} (1)_s (1)_{k-1-2t}} = 0
\end{equation}  
for all integers $i,t$, and the polynomial parameter $n$ and the rational function parameter $r$.
We denote the summand in \eqref{sum1} by $f_1(n,r,i,k,s,t)$ and observe by a simple calculation that  
\begin{equation}
\label{rec}  
f_1(n+6,r+12,i+3,k+6,s,t+3) = \frac{(1+r)_{12}}{64} f_1(n,r,i,k,s,t).
\end{equation} 
Let $g_1(n,r,i,t)$ denote the left-hand side of \eqref{sum1}. Summing \eqref{rec}, it follows that
\begin{equation}
\label{rec1} 
g_1(n+6,r+12,i+3,t+3) = \frac{(1+r)_{12}}{64} g_1(n,r,i,t).
\end{equation} 
Assume that we would have proven $g_1(n,0,i,t)=0$ for all integers $i,t$ and where $n$ is treated as a variable. We can then deduce from \eqref{rec1} that $g_1(n,0,i,t)=g_1(n,12,i,t)=g_1(n,24,i,t)=g_1(n,36,i,t)=\ldots=0$, and, since $g_1(n,r,i,t)$ is for fixed integers $n,i,t$ a rational function in $r$, this would prove $g_1(n,r,i,t)=0$ for all integers $n,i,t$. By the polynomiality in $n$, this would also imply that 
$g_1(n,r,i,t)=0$ is the zero polynomial.

In the remainder of this subsection we prove $g_1(n,0,i,t)=0$. By setting $r=0$ on the LHS of \eqref{sum1}, we obtain
\begin{multline*}
 \sum_{k,s} (-1)^{k+s} 2^{-k-s} \frac{(k-n)_s (0)_{k-1+2t} (-1/2-i+2k-n/2+s)_{2i-k-s}}{(1)_{2i-k-s} (1)_s (1)_{k-1-2t}} \\
= \sum_{k,s} (-1)^{1+s} 2^{-k-s} \frac{(k-n)_s  (-1/2-i+2k-n/2+s)_{2i-k-s}}{(1)_{2i-k-s} (1)_s (1)_{k-1-2t} (1)_{1-2t-k}}. 
\end{multline*} 
In order to prove that this sum evaluates to $0$, 
 we use Koutschan's Mathematica package HolonomicFunctions \cite{koutschan,christoph}, which is available at 
\url{https://www3.risc.jku.at/research/combinat/software/ergosum/RISC/HolonomicFunctions.html}.
Denote by $f_2(n,i,k,s,t)$ the summand.
Applying the command 
$$
{\tt FindCreativeTelescoping[f_2[n,i,k,s,t],\{S[k]-1,S[s]-1\},\{S[n],S[i],S[t]\}]}
$$ 
returns five polynomials $p_h$, $1 \le h \le 5$, for which the package provides the 
following identity
\begin{multline} 
\label{eq:rec via telescoping}
p_1(n,i,t) f_2(n,i+1,k,s,t) + p_2(n,i,t) f_2(n,i,k,s,t+1) + p_3(n,i,t) f_2(n,i,k,s,t) \\
+ \Delta_k  \left[ \frac{p_4(n,i,k,s,t)}{(1+2i-n)}\cdot  \frac{f_2(n,i,k,s,t)}{(2+2i-k-s)} \right] \\
+ \Delta_s \left[\frac{p_5(n,i,k,s,t)}{(1+2i-n)(n-k)}\cdot  \frac{f_2(n,i,k,s,t)}{(1+2i-k-s) (2+2i-k-s)} \right]=0, 
\end{multline}
where $\Delta_x$ denotes the forward difference, that is $\Delta_x d(x) = d(x+1)-d(x)$, 
and
\begin{align*}
p_1(n,i,t)&=4 (2 i+n-1) (2 i-2 t+1) \left(4 i^2+4 i n+32 i t+8 i+n^2-16 n t-8 n+8 t+3\right), \\
p_2(n,i,t)&=-64, \\
p_3(n,i,t)&=(2 i-n+4 t+1) \left(8 i^3+12 i^2 n+32 i^2 t+12 i^2+6 i n^2-32 i n t-12 i n-128 i t^2-48 i t-2
 i\right.\\
   &\left.+n^3-24 n^2 t-9 n^2+64 n t^2+72 n t+11 n+256 t^3+192 t^2+16 t-3\right),
\end{align*}
and $p_4,p_5$ are polynomials with 750 and 1168 monomials, respectively, which can be found in the accompanying Mathematica file {\tt companion.nb}. The file is available from the authors' webpages. It also checks 
the identity with Mathematica.

\begin{remark}
\label{remark_problem} 
\it 
The proofs of identities as \eqref{eq:rec via telescoping} are conceptually simple and work essentially as follows. We can use the two standard identities for Pochhammer symbols
\begin{equation}
\label{pochhammer_standard} 
(a)_n = a (a+1)_{n-1} \quad \text{and} \quad (a)_n = (a)_{n-1} (a+n-1) 
\end{equation}
to express quotients of ``shifts'' of $f(n,i,k,s,t)$ (such as for instances $f(n,i+1,k,s,t)$) in terms of 
$f(n,i,k,s,t)$ itself which is multiplied by an appropriate rational factor. This reduces the problem then to proving a certain rational function identity.

However, the application of \eqref{pochhammer_standard} can be problematic as it can cause divisions by zero. To give an example, observe that the second identity would imply $1=(1)_0 = (1)_{-1} 0 = \frac{0}{0}$ and it is then a matter of interpreting $\frac{0}{0}$ correctly. Before we argue that we can avoid divisions by zero, we present an example that shows how one can run into difficulties.

Consider $(n+1)_n$, then we can use the second identity twice to see that 
$
(n+1)_n = (n+1)_{n-2} (2n-1) (2n). 
$
On the other hand, 
$
(n)_{n-1} = n (n+1)_{n-2}  
$
by the first identity, so that we could be tempted to conclude that 
$
(n+1)_n = 2 (2n-1) n (n+1)_{n-2} = 2 (2n-1) (n)_{n-1}. 
$
However, $(n+1)_n = 2 (2n-1) (n)_{n-1}$ is simply not true for $n=0$ as $(1)_0=1$ and 
$(0)_{-1} = -1$. This is because there are problematic applications of  
\eqref{pochhammer_standard}: $(n+1)_n = (n+1)_{n-1} (2n)$ and $(n)_{n-1} = n (n+1)_{n-2}$ 
lead for $n=0$ to indefinite expressions of the form $\frac{0}{0}$ on the right-hand side.  

We can circumvent such problems by working with the limit
$(a)_m = \lim_{\epsilon \to 0} (a+\epsilon)_m$. In our example, this would lead to the following.
\begin{multline*}
(n+1)_n  = \lim_{\epsilon \to 0} (n+\epsilon)_{n-1} \frac{(n+1+\epsilon)_n}{(n+\epsilon)_{n-1}}  = \lim_{\epsilon \to 0} (n+\epsilon)_{n-1} \frac{(2n-1+\epsilon)(2n+\epsilon)}{n+\epsilon}  \\
= \lim_{\epsilon \to 0} (n+\epsilon)_{n-1} \frac{(2n-1+\epsilon)(2n+\epsilon)}{n+\epsilon}  = 
\lim_{\epsilon \to 0} (n+\epsilon)_{n-1} (2n-1+\epsilon) \frac{2n+\epsilon}{n+\epsilon} = 
 (n)_{n-1} (2n-1) \lim_{\epsilon \to 0} \frac{2n+\epsilon}{n+\epsilon}.
\end{multline*} 
Now observe that 
$$
\lim_{\epsilon \to 0} \frac{2n+\epsilon}{n+\epsilon} = \begin{cases} 1 & n=0 \\ 2 & n \not= 0 \end{cases}, 
$$
so that $(n+1)_n = (1+[n \not= 0]) (2n-1) (n)_{n-1}$, using Iverson bracket, which is defined as 
$[\text{statement}]=1$ if the statement is true and $[\text{statement}]=0$ otherwise.

Suppose we need to prove an identity such as \eqref{eq:rec via telescoping} and we want to relate shifts of the summand to the non-shifted summand via the multiplication of a rational function. For this purpose, we obtain for individual Pochhammer symbols $(a)_m$ in the summand expressions of the form
$$(a+s)_{m+t} = \lim_{\epsilon \to 0} (a+\epsilon+s)_{m+t} = \lim_{\epsilon \to 0}  (a+\epsilon )_m \frac{(a+\epsilon+s)_{m+t}}{(a+\epsilon)_m} =  \lim_{\epsilon \to 0} (a+\epsilon)_m \frac{(a+\epsilon+m)_{s+t}}{(a+\epsilon)_{s}},$$ 
where $s,t$ are fixed integers of small absolute value (they come from the shifts of the parameters). Note that this also takes care of Pochhammer symbols in the denominator 
using $1/(a)_{-m}=(a-m)_{m}$.

Now one may run into difficulties in the sense that we cannot easily remove the limit, since some linear factors $a+j$ of the denominator of $\frac{(a+m)_{s+t}}{(a)_{s}}$ vanish for some choices of parameters.
However, since we assume that the summand is well-defined for all choices of parameters, these are always removable singularities. 
More concretely, assume that the summand is 
$$
f(a_1,\ldots,a_n;m_1,\ldots,m_n) =  \prod_{i=1}^n (a_i)_{m_i},
$$
then we would like to replace the shifted summand
$f(a_1+s_1,\ldots,a_n+s_n;m_1+t_1,\ldots,m_n+t_n)$ by 
\begin{equation}
\label{rel}  
\lim_{\epsilon \to 0}  \prod_{i=1}^n (a_i+\epsilon)_{m_i} \frac{(a_i+ \epsilon + m_i)_{s_i+t_i}}{(a_i + \epsilon )_{s_i}} = \lim_{\epsilon \to 0} f(a_1+\epsilon,\ldots,a_n+\epsilon;m_1,\ldots,m_n) \prod_{i=1}^n \frac{(a_i+\epsilon + m_i)_{s_i+t_i}}{(a_i + \epsilon)_{s_i}}.
\end{equation} 
In the following, we assume that $s_i \ge 0$ and $s_i+t_i \ge 0$ for all $i$, but the argument is analogous for the other cases.

Then a singularity $a_i+j$ in the denominator of the rational function on the right-hand side of \eqref{rel} can clearly only be removable for the following two reasons: either it is removable due to $f(a_1+\epsilon,\ldots,a_n+\epsilon;m_1,\ldots,m_n)$ or 
due to $\prod_{i=1}^n (a_i+\epsilon + m_i)_{s_i+t_i}$. The first case is dealt with in (1), while the second case is dealt with in (2).

\begin{enumerate} 
\item The factor $a_i+j$ ``cancels'' with $f(a_1,\ldots,a_n;m_1,\ldots,m_n)$, that is if $a_i+j=0$ then also $f(a_1,\ldots,a_n;m_1,\ldots,m_n)=0$. Then one may simply multiply the rational function identity with all such linear factors to obtain a rational function identity that is valid for all choices of parameters, and this multiplication is justified as the linear factor appears as part of the Pochhammer symbols in the summand.
\item The second option is that the factor $a_i+j$ cancels with a factor in some $(a_k+m_k)_{s_k+t_k}$.
This is precisely the situation where Iverson notation may come into the play. Concretely, it may happen that there is a $k \in \{1,2,\ldots,n\}$
and an integer $l$ such that
\begin{equation}
\label{implies}  
a_i+j=0 \Rightarrow a_k+m_k+l=0.
\end{equation} 

In the current paper, the $a_i$'s and the $m_i$'s are typically affine combinations of the parameters over the integers.
Under this circumstance, we have \eqref{implies} only if there is a rational number $q$ with $q(a_k+m_k+l)=a_i+j$. It can be checked easily that this second case never happens in the applications in this paper.
\end{enumerate} 

Also note the following subtlety: By our assumption that the $a_i$'s and the $m_i$'s are affine combinations of the parameters, $a_i+j=0$ is an affine hyperplane in these parameters, and, in principle, it may happen that for different choices of parameters we are either in Case (1) or in Case (2). In addition, in Case (2) we might even have to consider different terms $a_k+m_k+l$ for different choices of parameters. However, this is actually not possible because the terms in the numerators that are responsible for the fact that the zeros are removable also induce affine hyperplanes and an affine hyperplane $H$ is never contained in the union of a finite number of affine hyperplanes that are  different from $H$.

Finally note that all these considerations can easily be extended to summands of the form 
$\prod_{i=1}^k z_i^{l_i} \prod_{i=1}^n (a_i)_{m_i}$, where the $z_i$ are integers and the $l_i$ are also affine linear combinations of the parameters over $\mathbb{Z}$. 
\end{remark}
 
Now we sum \eqref{eq:rec via telescoping} over all $k,s,t \in \mathbb{Z}$. Let us repeat that, since we treat $n$ as a variable, the singularities caused by the factors $1+2i-n$ and $n-k$ in the numerator do not lead to any problems. As for $1+2i-k-s$ and $2+2i-k-s$, note that $f_2(n,i,k,s,t)$ 
vanishes if $1+2i-k-s=0$ and $2+2i-k-s=0$, which is due to the factor $1/(1)_{2i-k-s}$, therefore also these singularities are removable.

As a consequence of  the telescoping effect on the forward difference terms in \eqref{eq:rec via telescoping}, we can conclude 
\begin{equation}
\label{rec4} 
p_1(n,i,t) g_1(n,0,i+1,t) + p_2(n,i,t) g_1(n,0,i,t+1) + p_3(n,i,t) g_1(n,0,i,t)=0.
\end{equation} 
Now we can prove $g_1(n,0,i,t)=0$ by induction with respect to $i-t$. Since $f_2(n,i,k,s,t)$ is only non-zero for $1+2t \leq k \leq 1-2t$ and $0 \leq s \leq 2i-k$, we can conclude that $g_1(n,0,i,t)=0$ for $i \le t$. We assume by the 
induction hypothesis that $g_1(n,0,i,t+1)=0$ and $g_1(n,0,i,t)=0$, so that, by \eqref{rec4}, $p_1(n,i,t) g_1(n,0,i+1,t)=0$. 
Now we need to argue that $p_1(n,i,t)$ is non-zero. For integers $i,t$, all three factors in $p_1(n,i,t)$ are non-zero polynomials in $n$, since their constant terms are odd numbers.

\subsection{An auxiliary double sum} 
\label{sec:auxiliary double sum}
We let $f_3(n,r,i,j,k,s,t)$ denote the summand of the triple sum in \eqref{triple}, i.e.,
\begin{multline*} 
f_3(n,r,i,j,k,s,t) = (-1)^{1+k+s+t} 2^{2-k-s} \\ 
\times \frac{(k-n)_s (-r)_{k-1+2t} (-1/2-i+2k-n/2-r/2+s)_{2i-k-s} (2-2j+r)_{j-1-2t} (j-t)_{j-t-1}}{(1)_{2i-k-s} (1)_s (1)_{k-1-2t} (1)_t}.
\end{multline*} 
The aim of this subsection is to prove
\begin{equation}
\label{eq:aux sum}
\sum_{s,t} f_3(n,r,i,j,j+1,s,t) = 0,
\end{equation}
where $n$ is treated as a polynomial parameter, $r$ is treated as a rational function parameter, and $i$ is a non-negative integer and $j \ge 2$ is an integer. We show that, in fact, the inner sum of \eqref{eq:aux sum} over $t$ evaluates to zero. Neglecting all terms that do not depend on $t$, we need to show 
$$
\sum_{t} (-1)^t \frac{(-r)_{j+2t} (2-2j+r)_{j-1-2t} (j-t)_{j-t-1}}{(1)_{j-2t} (1)_t}=0.
$$
Rewritten in hypergeometric series notation, this is 
\begin{equation}
\label{eq:aux in hyp}
\frac{(j)_{j-1} (-r)_j (2-2j+r)_{j-1}}{(1)_j} \pFq{2}{1}{-\frac{j}{2},\frac{1}{2}-\frac{j}{2}}{\frac{3}{2}-j}{1}=0.
\end{equation}
Now we use the Chu-Vandermonde summation \cite[(1.7.7); Appendix (III.4)]{slater},
\begin{equation}
\label{eq:S2101}
\pFq{2}{1}{a,-n}{c}{1} = \frac{(c-a)_n}{(c)_n},
\end{equation}
where $n$ is a non-negative integer, and in our application we need to distinguish between $j$ is even and $j$ is odd. Namely, for $j=2a$ and $j=2a+1$, respectively, the term $(c-a)_n$ in the corresponding summation becomes $(1-a)_a$, which is equal to $0$ and hence implies \eqref{eq:aux sum}.
 
\subsection{Proof of \eqref{triple}}
\label{sec:proof of triple sum}

Recall that we denote the summand in \eqref{triple} by $f_3(n,r,i,j,k,s,t)$. We will prove the following recursion for the summand, which is true unless $j=t$.
\begin{multline} 
\label{rec3}
(3j-r-2)(r+1)_4 f_3(n,r,i,j,k,s,t) = \\
2 (2j+1) (2-2j+r)_2 f_3(n+2,r+4,i+1,j+1,k+2,s,t+1) \\
+ (j-r-3) f_3(n+2,r+4,i+1,j+2,k+2,s,t+1) 
\end{multline} 
In our application, this exception is not a problem, since we can assume $t < j$ as $t \le (k-1)/2$ 
(otherwise the summand vanishes due to the factor $1/(1)_{k-1-2t}$) 
 and we have the condition $k \le j$ from the range of our sum \eqref{triple}. Indeed,
\begin{multline*}
(3j-r-2)(r+1)_4 f_3(n,r,i,j,k,s,t) =
(3j-r-2)(r+1)_4 (-1)^{1+k+s+t} 2^{2-k-s} \\ 
\times \frac{(k-n)_s (-r)_{k-1+2t} (-1/2-i+2k-n/2-r/2+s)_{2i-k-s} (2-2j+r)_{j-1-2t} (j-t)_{j-t-1}}{(1)_{2i-k-s} (1)_s (1)_{k-1-2t} (1)_t}\\
= - 4\frac{(3j-r-2)(t+1)(2-2j+r)_2}{(r+1-j-2t)}f_3(n+2,r+4,i+1,j+1,k+2,s,t+1),
\end{multline*}
and 
\begin{multline*}
(j-r-3) f_3(n+2,r+4,i+1,j+2,k+2,s,t+1) 
=  (j-r-3)(-1)^{k+s+t} 2^{-k-s} \\ \times 
\frac{(k-n)_s (-r-4)_{k+3+2t} (-1/2-i+2k-n/2-r/2+s)_{2i-k-s} (2-2j+r)_{j-1-2t} (j-t+1)_{j-t}}{(1)_{2i-k-s} (1)_s (1)_{k-1-2t} (1)_{t+1}} \\
= \frac{(j-r-3)(2-2j+r)_2(2j-2t)(2j-2t-1)}{(r+1-j-2t)(j-t)} f_3(n+2,r+4,i+1,j+1,k+2,s,t+1)\\
= \frac{2(j-r-3)(2-2j+r)_2(2j-2t-1)}{(r+1-j-2t)} f_3(n+2,r+4,i+1,j+1,k+2,s,t+1),
\end{multline*} 
where we need to assume $j\not= t$ in the last step. Now
\begin{multline*} 
2 (2j+1) (2-2j+r)_2 f_3(n+2,r+4,i+1,j+1,k+2,s,t+1) \\
+ (j-r-3) f_3(n+2,r+4,i+1,j+2,k+2,s,t+1) \\ 
= -4\frac{(3j-r-2)(t+1)(2-2j+r)_2}{(r+1-j-2t)}f_3(n+2,r+4,i+1,j+1,k+2,s,t+1) \\
=(3j-r-2)(r+1)_4 f_3(n,r,i,j,k,s,t).
\end{multline*}

We sum the recursion \eqref{rec3} over all $s,t \in \mathbb{Z}$. Letting 
$$
g_3(n,r,i,j,k) = \sum_{s,t} f_3(n,r,i,j,k,s,t),
$$
we have 
\begin{multline*} 
(3j-r-2)(r+1)_4 g_3(n,r,i,j,k) = \\
2 (2j+1) (2j-r-3)_2 g_3(n+2,r+4,i+1,j+1,k+2) 
+ (j-r-3) g_3(n+2,r+4,i+1,j+2,k+2).
\end{multline*} 
Now we sum over all $k \le j$ and let 
$$
h_3(n,r,i,j) = \sum_{k \le j} g_3(n,r,i,j,k).
$$
We obtain 
\begin{multline*} 
(3j-r-2) (r+1)_4 h_3(n,r,i,j) = \\
2 (2j+1) (2j-r-3)_2 \Big[ h_3(n+2,r+4,i+1,j+1) - g_3(n+2,r+4,i+1,j+1,j+2) \Big] \\
+(j-r-3) h_3(n+2,r+4,i+1,j+2).
\end{multline*}
Assuming $j \ge 1$, the auxiliary double sum evaluation  \eqref{eq:aux sum} shows that $g_3(n+2,r+4,i+1,j+1,j+2)=0$ so that we have 
\begin{multline*} 
(3j-r-2) (r+1)_4 h_3(n,r,i,j) = \\
2 (2j+1) (2j-r-3)_2 h_3(n+2,r+4,i+1,j+1) 
+(j-r-3) h_3(n+2,r+4,i+1,j+2).
\end{multline*}

Now we compute $h_3(n,r,i,j)$ for $i<j$ by induction with respect to $i$.  The base case $i=1$ follows immediately from Section~\ref{sec:j large} since $2i=2 \leq j$. In the induction step, we use a second induction with respect to $j-i$, where we go from large $j-i$ to small. The base case $2i\leq j$ follows again from  Section~\ref{sec:j large}. For the induction step, we can assume $h_3(n,r,i,j)=0$ (by the induction hypothesis of the first induction) and $h_3(n+2,r+4,i+1,j+2)=0$ (by the induction hypotheses of the second induction). Hence we have $2 (2j+1) (2j-r-3)_2 h_3(n+2,r+4,i+1,j+1) =0$ which implies $h_3(n+2,r+4,i+1,j+1) =0$ since we treat $r$ as a variable.

For $i=j$, the above recursion simplifies to
\[
(3i-r-2) (r+1)_4 h_3(n,r,i,i) = 
2 (2i+1) (2i-r-3)_2 h_3(n+2,r+4,i+1,i+1).
\]
Again we use induction with respect to $i$. We can check the base case $i=1$ easily since $f_3(n,r,1,1,k,s,t)=0$ unless $s,t$ are non-negative and satisfy $k+s\leq 2$ and $2t+1\leq k$. For the induction step, we use the above recursion.

This concludes the proof of Lemma~\ref{prop:w1 matrix product}, and hence of Theorem~\ref{1}.

\section{Proof of Theorem~\ref{3}}
\label{sec:3}

\subsection{LU-decomposition} 

On specializing $w=0$ in \eqref{odddet}, we obtain the following.
$$
6^n \det_{1 \le i, j \le n} \left( \sum_{p,q} 0^{n-j-q} (-1)^j \binom{j-1}{p} \binom{n-j}{q}  
\binom{p-q-l+i-2}{2i-1} \right), 
$$
Since $0^{n-j-q}=0$ unless $q=n-j$, the above simplifies to 
$$
6^n \det_{1 \le i, j \le n} \left( \sum_{p} (-1)^j \binom{j-1}{p}  
\binom{p-n+j-l+i-2}{2i-1} \right).
$$
We define 
$$
b_{i,j} = \sum_{p} (-1)^j \binom{j-1}{p}  \binom{p-n+j-l+i-2}{2i-1}.
$$

Also here we were able to guess the LU-decomposition.  Namely, let 
\[
y_{i,j} = \sum_{t} 2(-1)^{t+j}3^{1-j}\frac{(j)_j}{(1)_{j-t}(1)_{t-i}(3-i+2n+2l-2t)_{i}}.
\]

The remainder of this section is devoted to the proof of the following lemma.

\begin{lemma}
\label{prop:w0 matrix product}
Let $n$ be a positive integer, then $(b_{i,j})_{1 \le i, j \le n} \cdot (y_{i,j})_{1 \le i,j \le n}$ is a lower triangular 
matrix with $1$'s on the main diagonal, that is
\begin{equation}
\label{eq:w0 matrix identity}
\sum_{k=1}^n b_{i,k}  y_{k,j} = \begin{cases}
1 \qquad &i=j,\\
0 & i<j.
\end{cases}
\end{equation}
\end{lemma}

The lemma implies
$$
\det_{1 \le i, j \le n} \left( \sum_{p} (-1)^j \binom{j-1}{p}  
\binom{p-n+j-l+i-2}{2i-1} \right) = \left[ \det_{1 \le i, j \le n} (y_{i,j}) \right]^{-1}
$$
and it is easy to see that  $6^n \left[ \det_{1 \le i, j \le n} (y_{i,j}) \right]^{-1}$ is the expression in Theorem~\ref{3} when setting $m=2l+1$ there.

Observe that the left-hand side in \eqref{eq:w0 matrix identity} is a rational function in $l,n$, in fact, it is a rational function in $l+n$, and thus it suffices to consider the case $l=0$, that is, we need to show
\begin{equation}
\label{eq:w0 equation}
\sum_{k \ge 1,s,t}  \frac{(-1)^{k+s+t} (1-k)_s (-i+k-n+s)_{2i-1}}{(1)_s (1)_{j-t} (1)_{t-k} (3-k+2n-2t)_k}
= \begin{cases} \frac{1}{2} (-1)^j 3^{j-1} (1)_{j-1} & 1 \le i=j \\ 0 & 1\le i < j \end{cases},
\end{equation}
where we have multiplied both sides by $\frac{(-1)^j 3^{j-1}(1)_{2i-1}}{2(j)_j} $ and have replaced the summation parameter $p$ by $s$.

\subsection{The case $i=j$}
Denote by $f_4(n,i,j,k,s,t)$ the summand on the LHS of \eqref{eq:w0 equation}.
We use again Koutschan's Mathematica package HolonomicFunctions.
Applying the command 
$$
{\tt FindCreativeTelescoping[f_4[n,i,i,k,s,t],\{S[k]-1,S[s]-1,S[t]-1\},\{S[i],S[n]\}]}
$$ 
returns seven polynomials $p_h$, $1 \le h \le 7$, where the first four do not depend 
on $k,s,t$, for which the package provides the following identity
\begin{multline}
\label{eq:rec for ii}
p_1 f_4(n+1,i+1,i+1,k,s,t)+p_2 f_4(n+1,i,i,k,s,t)+p_3 f_4(n,i+1,i+1,k,s,t)+p_4 f_4(n,i,i,k,s,t)\\
+\Delta_k \left[\frac{(k-s-1) \, p_5 }{(1 + i - t)} \cdot \frac{f_4(n,i,i,k,s,t)}{2q(-2 + i + k - n + s) }\right]
+\Delta_s \left[\frac{s \cdot p_6 }{(1 + i - t)} \cdot \frac{f_4(n,i,i,k,s,t)}{2q(-2 + i + k - n + s) }\right] \\
+\Delta_t \left[\frac{(k-t) (k-2n+2t-3) \, p_7}{(1 + i - t)} \cdot \frac{f_4(n,i,i,k,s,t)}{q}\right] =0,
\end{multline}
where $q=2i(-9 + 8 i + i^2 - 8 n - 2 i n + n^2)(3 + 2 n - 2 t)(2 + n - t)$ and
\begin{align*}
p_1 &= -6 i + 14 i^2 + 3 n - 24 i n + 7 n^2,\\
p_2 &= 2 i (-1 + 7 i - 3 n) (i - n),\\
p_3 &= -(-1 + 3 i - 2 n) (3 i - 2 n),\\
p_4 &= i (-7 i + i^2 + n - 16 i n + 3 n^2),
\end{align*}
and $p_5,p_6,p_7$ are polynomials with 1031, 1031, and 312 monomials, respectively, which can be found in the accompanying Mathematica file  {\tt companion.nb}. 

Denote by $g_4(n,i,j) = \sum_{k \geq 1}\sum_{s,t}f_4(n,i,j,k,s,t)$. Since $n$ is a polynomial variable, the singularities caused by terms that involve $n$ are removable singularities. As for $(1+i-t)$, note that $f_4(n,i,i,k,s,t)$ vanishes if $1+i-t=0$, which is due to the factor $1/(1)_{i-t}$.   

As a consequence of the telescoping in $s,t$ and $k$ the above recurrence implies
\begin{multline*}
p_1 g_4(n+1,i+1,i+1)
+p_2 g_4(n+1,i,i)
+p_3 g_4(n,i+1,i+1)
+p_4 g_4(n,i,i)\\
= \left[ \sum_{s,t}\frac{\left. (-s) p_5\right|_{k=1}}{(1 + i - t)} \cdot \frac{f_4(n,i,i,1,s,t)}{2q(-1 + i - n + s) }\right].
\end{multline*}
By definition $g_4(n,i,i,1,s,t)=0$ unless $s=0$. Hence the above recursion simplifies to
\begin{equation}
\label{eq:rec for ii simp}
p_1 g_4(n+1,i+1,i+1) +p_2 g_4(n+1,i,i) +p_3 g_4(n,i+1,i+1) +p_4 g_4(n,i,i)=0.
\end{equation}
We complete the proof in the $i=j$ case by showing $g_4(n,i,i)= \frac{1}{2}(-1)^i 3^{i-1}(i-1)!$ by induction with respect to $i$. The base case is $i=1$. It is immediate that $f_4(n,1,1,k,s,t)=0$ unless $t=k=1$ and $s=0$. Hence $g_4(n,1,1)=f_4(n,1,1,1,0,1)=-\frac{1}{2}$.

In the induction step, we first show that $g_4(n,i+1,i+1)$ is independent of $n$. By definition, $g_4$ is a rational function in $n$, that is,  $g_4(n,i+1,i+1)=a(n)/b(n)$, where $a(n),b(n)$ can assumed to be polynomials in $n$ without a common factor. By \eqref{eq:rec for ii simp} and the induction hypotheses we obtain
\begin{equation}
\label{eq:help ii}
p_1 a(n+1)b(n)+p_3 a(n)b(n+1) = \frac{1}{2}(-3)^i i! b(n)b(n+1) (3i-5i^2-n+12i n-3n^2).
\end{equation}
Hence, $b(n)$ divides $p_3 b(n+1)=-(2n-3i)(2n-3i+1)b(n+1)$. 
It is not difficult to see that the sum in the definition of $g_4$ ranges over all $1 \leq s+1 \leq k \leq t \leq i$. Hence the possible factors of $b(n)$ are $(2n-3i+3),(2n-3i+4),\ldots,(2n)$, which implies that $b(n)$ and $p_3$ can not have a common factor. Therefore $b(n)$ divides $b(n+1)$, and this implies $b(n)=1$. By comparing the degrees on both sides of \eqref{eq:help ii} and taking the leading coefficients of $p_1$ and $p_3$ as polynomials in $n$ into account, we see that $a(n)$ must be a polynomial in $n$ of degree $0$, i.e., $a(n)=a(n+1)$. The assertion now follows immediately from \eqref{eq:rec for ii simp}, where we also use that  $p_1+p_3$ is non-zero since the constant term as polynomial in $n$ is positive for $i \ge 1$.

\subsection{The case $i<j$} 
As before, we use the command
$$
{\tt FindCreativeTelescoping[f_4[n,i,j,k,s,t],\{S[k]-1,S[s]-1,S[t]-1\},\{S[i],S[j],S[n]\}]}
$$ 
and obtain seven polynomials $p_h(n,i,j,k,s,t)$, $1 \le h \le 7$, where the first four do not 
depend on $k,s,t$, 
for which the package provides the following identity
\begin{multline}
\label{eq:rec for ij}
p_1 f_4(n+1,i,j+1,k,s,t)
+p_2 f_4(n+1,i,j,k,s,t)
+p_3 f_4(n,i,j+1,k,s,t)
+p_4 f_4(n,i,j,k,s,t) \\
+\Delta_k \left[ \frac{p_5}{(-1-j+t)} \cdot \frac{f_4(n,i,j,k,s,t)}{2 (-2+i+k-n+s)(3+2n-2t)(2+n-t) } \right]\\
+\Delta_s \left[ \frac{s \cdot p_6}{2(k-s)(-1-j+t)}\cdot  \frac{f_4(n,i,j,k,s,t)}{(3+2n-2t)} \right]\\
+\Delta_t \left[ \frac{(k-t)(-3+k-2n+2t) \, p_7}{(k-s)(-1-j+t)}\cdot  \frac{f_4(n,i,j,k,s,t)}{2(-i+k-n+s)(3+2n-2t)(2+n-t)} \right]
=0,
\end{multline}
where
\begin{align*}
p_1&= 3 i - i^2 - 10 i j + 2 i^2 j + 6 i j^2 - n + 9 i n - i^2 n + 6 j n - 
 12 i j n - 6 j^2 n - 2 n^2 + 6 i n^2 + 6 j n^2 - n^3,\\
p_2&= -2 (-2 i + i^2 + 3 i j - 3 i n - 3 j n),\\
p_3&= (-1 + 3 j - 2 n) (3 j - 2 n) (i - 2 j + n),\\
p_4&= 2 i - i^2 - 6 j - 9 i j + 18 j^2 + 3 i n - 15 j n.
\end{align*}
and $p_5,p_6,p_7$ are polynomials with 484, 124, and 173 monomials, respectively, which can be found in the accompanying Mathematica file {\tt companion.nb}. 

As before, the above recursion implies
\begin{multline}
\label{eq:rec for ij 2}
p_1 g_4(n+1,i,j+1)
+p_2 g_4(n+1,i,j)
+p_3 g_4(n,i,j+1)
+p_4 g_4(n,i,j)\\
= \left[ \sum_{s,t}\frac{\left.p_5\right|_{k=1}}{(-1 -j +t)} \cdot \frac{f_4(n,i,j,1,s,t)}{2 (-1+i-n+s)(3+2n-2t)(2+n-t) }\right].
\end{multline}
It is immediate that $f_4(n,i,j,1,s,t)=0$ unless $s=0$.
The right hand side can therefore be rewritten as
\[
\sum_{r=1}^3 \frac{c_r(n,i,j)(-i+1-n)_{2i-1}}{4(n+1-i)} \sum_{t}\frac{(-1)^{t+1}}{(1)_{j-t+1}(1)_{t-r}(1+n-t)(2+n-t)},
\]
where the $c_r(n,i,j)$ are defined through $\frac{\left.p_5\right|_{k=1}}{(3+2n-2t)}=c_1(n,i,j)+(t-1)c_2(n,i,j)+(t-1)(t-2)c_3(n,i,j)$. 
By rewriting the sum over $t$ in hypergeometric notation and using the Chu-Vandermonde summation, see \eqref{eq:S2101}, we obtain
\begin{multline*}
\sum_{t}\frac{(-1)^{t+1}}{(1)_{j-t+1}(1)_{t-r}(1+n-t)(2+n-t)} \\
= \frac{(-1)^{r+1}}{(1+n-r)(2+n-r)(1)_{j+1-r}}\pFq{2}{1}{-2-n-r,-1-j+r}{-n+r}{1}
= \frac{(-1)^{r+1}(2+j-r)}{(-n+r-2)_{j+3-r}},
\end{multline*}
where we assume $j+1-r \geq 0$ which is always true since $j\geq 2$ and $r \leq 3$.
Using the explicit values for the $c_i(n,i,t)$
\begin{align*}
c_1(n,i,j)&= -2 (i + n) (-2 - i + i^2 + 6 j - 5 n - 2 i n - 3 n^2),\\
c_2(n,i,j)&= -2 (i + n) (-i + i^2 - i j + i^2 j + 3 n + 2 i n - 2 i^2 n - 13 j n + 
   9 n^2 + 4 i n^2 - j n^2 + 6 n^3),\\
c_3(n,i,j)&= -2 (1 + n) (i + n) (i j - i^2 j - 2 n - i n + i^2 n + 7 j n - 5 n^2 - 
   2 i n^2 + j n^2 - 3 n^3),
\end{align*}
we obtain
\[
\sum_{r=1}^3 \frac{(-1)^{r+1}c_r(n,i,j)(-i+1-n)_{2i-1}(2+j-r)}{4(n+1-i)(-n+r-2)_{j+3-r}}=0.
\]
Hence \eqref{eq:rec for ij 2} simplifies to
\begin{equation}
\label{eq:rec for ij simp}
p_1 g_4(n+1,i,j+1)
+p_2 g_4(n+1,i,j)
+p_3 g_4(n,i,j+1)
+p_4 g_4(n,i,j)=0
\end{equation}
Finally we prove by induction with respect to $j$ that $g_4(n,i,j)=0$ if $i<j$.

As for the base case $j=i+1$ of the induction, note that we obtain 
$$
p_1 g_4(n+1,i,i+1)=-p_3g_4(n,i,i+1)
$$
when setting $i=j$ in \eqref{eq:rec for ij simp} as $p_2+p_4=0$ under this assumption. We will work with the analogous recursion in the induction step for $j>i+1$, and so we can defer the further treatment also of the case $j=i+1$ to the induction step.

 For the induction step, we assume that $g_4(n,i,j)=0$, then the above recursion simplifies to 
\begin{equation}
\label{eq:rec for ij simp2}
p_1 g_4(n+1,i,j+1)=-p_3g_4(n,i,j+1).
\end{equation}
For fixed $i$ and $j$, it follows from the definition that $g_4(n,i,j+1)$ is a rational function in $n$, that is
\[
g_4(n,i,j+1)=\frac{a(n)}{b(n)},
\]
for some polynomials $a(n),b(n)$. We can assume that $b(n)$ is monic and denote the leading coefficient of $a(n)$ by $\alpha$. By multiplying \eqref{eq:rec for ij simp2} by $b(n)b(n+1)$ we obtain
\[
p_1 a(n+1)b(n) = - p_3 a(n) b(n+1).
\]
Now the leading coefficient on the left-hand side is $-\alpha$ while on the right-hand it is $-4\alpha$, hence $a(n)$ must be the zero polynomial and we have $g_4(n,i,j+1)=0$ which proves the claim.

\section{Acknowledgements} 

We thank Christoph Koutschan and Carsten Schneider for useful discussions. In particular, Carsten Schneider also demonstrated to us how \eqref{eq:w0 equation} can be shown with his Mathematica Package  Sigma \cite{schneider,carsten}, which is available at \url{https://www3.risc.jku.at/research/combinat/software/Sigma/}.

\end{document}